%% file: m3-I-2.tex
\input m3-macs

\pageno=19

\tinfo I.2.19-29

\SetTFLinebox{\gtp }
\SetFLinebox{\gtv3 }
\SetHLinebox{\issn}

\H 2. $p$-primary part of the Milnor $K$-groups and \\ 
Galois cohomologies of fields of characteristic $p$ 

Oleg Izhboldin

\SetAuthorHead{O. Izhboldin}
\SetTitleHead{Part I. Section 2. $K$-groups  and  
Galois cohomologies of fields of characteristic $p$ \qquad\qquad}

\HH 2.0. Introduction

Let $F$ be a field and $F^{\sep}$ be the separable closure of $F$.
Let $F^{\ab}$ be the maximal abelian extension of $F$.
Clearly the Galois group $G^{\ab}=\Gal(F^{\ab}/F)$
is canonically isomorphic to the quotient of the
absolute Galois group $G=\Gal(F^{\sep}/F)$ modulo
the closure of its commutant.
By  Pontryagin duality, 
a description of $G^{\ab}$ is equivalent to a description of 
$$
\Hom_{\,\text{cont}}
\bigl(G^{\ab}, \Bbb Z/m\bigr) =
\Hom_{\,\text{cont}}\bigl( G, \Bbb Z/m\bigr) 
=H^{1}(F, \Bbb Z/m).
$$
where $m$ runs over all positive integers.
Clearly, it suffices to consider the case where $m$ is  
a power of a prime, say $m=p^i$.
The main cohomological tool to
compute the group $H^{1}(F, \Bbb Z/m)$   is a pairing 
$$
(\,\,,\,\,)_{m}\colon H^{1}(F, \Bbb Z/m)\otimes
K_{n}(F)/m\rightarrow H_{m}^{n+1}(F)
$$
where the right hand side is a certain cohomological
group discussed below.

Here $K_n(F)$ for a field $F$ is the $n$th Milnor $K$-group 
$K_{n}(F)=K_{n}^{M}(F)$ defined as 
$$(F^{*})^{\otimes n}/J$$ 
where $J$ is the subgroup generated by the elements of the form 
$a_{1} \otimes ... \otimes a_{n}$ 
such that $a_{i}+a_{j}=1$ for some $i \neq j$. 
We denote by $\{a_{1},...,a_{n}\}$ the class of 
$a_{1} \otimes ... \otimes a_{n}$. 
Namely, $K_{n}(F)$ is the abelian group defined by the following 

\noindent generators: symbols $\{a_{1},...,a_{n}\}$ with
$a_{1}$,...,$a_{n}$ 
$\in F^{*}$ 

\noindent and relations: 
$$
\aligned
&\{a_{1},...,a_{i}a_{i}',...a_{n}\}   
= \{a_{1},...,a_{i},...a_{n}\} + 
\{a_{1},...,a_{i}',...a_{n}\}\\ 
&\{a_{1},...,a_{n}\}=0\quad\text{ if $a_{i}+a_{j}=1$ for some $i$ and $j$ with 
$i \neq j$.} \endaligned $$ 
We write the group law additively.

\bigskip

Consider the following  
example (definitions of the groups will be given later).

\eg Example

Let $F$ be a field and let $p$ be a prime integer.
Assume that there is an integer $n$ with the following properties:
\Roster
\Item{(i)} 
the group $H_p^{n+1}(F)$ is isomorphic to $\Bbb Z/p$,
\Item{(ii)}
the pairing 
$$
(\,\,,\,\,)_{p}\colon H^{1}(F, \Bbb Z/p)\otimes
K_{n}(F)/p\rightarrow H_{p}^{n+1}(F)\simeq \Bbb Z/p
$$ 
is non-degenerate in a certain sense.
\endRoster
\noindent Then the $\Bbb Z/p$-linear space $H^1(F,\Bbb Z/p)$ is obviously dual to the $\Bbb Z/p$-linear space $K_{n}(F)/p$. 
On the other hand, $H^1(F,\Bbb Z/p)$ is dual to the $\Bbb Z/p$-space
$G^{\ab}/(G^{\ab})^p$. 
Therefore there is an isomorphism 
$$
\Psi_{F,p}\colon K_n(F)/p\simeq G^{\ab}/(G^{\ab})^p.
$$
\endeg

It turns out that this example 
can be applied to computations of the group $G^{\ab}/(G^{\ab})^p$ for 
multidimensional local fields. Moreover, it is possible to show that
the homomorphism $\Psi_{F,p}$ can be naturally extended to a homomorphism
$\Psi_F\colon K_n(F)\to G^{\ab}$ (the so called reciprocity map).
Since $G^{\ab}$ is a profinite group, it follows that the homomorphism
$\Psi_F\colon K_n(F)\to G^{\ab}$ factors through the homomorphism
$K_n(F)/DK_n(F)\to G^{\ab}$ where the group $DK_n(F)$ consists of
all divisible elements:  $$DK_n(F):=\cap_{m\ge 1} mK_n(F).$$
This observation makes natural the following notation:

\df Definition {{\rm(cf.~section 6 of Part~I)}}

For a field $F$ and integer $n\ge0$ set
$$
K_n^{t}(F):=K_n(F)/DK_n(F),
$$ 
where $DK_n(F):=\bigcap_{m\ge 1} mK_n(F)$.
\enddf

The group $K_n^{t}(F)$ for a higher local field $F$ endowed with a certain topology 
(cf. section 6 of this part of the volume) 
is called a  topological Milnor $K$-group
$K^{\tpp}(F)$ of $F$. 

The example shows that computing  the group $G^{\ab}$ is closely 
related to computing the groups $K_n(F)$, $K_n^{t}(F)$, 
and $H^{n+1}_m(F)$.
The main purpose of this section is to explain some basic properties of
these groups and discuss several classical conjectures.
Among the problems, we point out the following:

\bitem
discuss $p$-torsion and cotorsion  of the groups 
$K_n(F)$ and $K_n^{t}(F)$,
\bitem
study  an analogue of Satz 90  for the
groups $K_n(F)$ and $K_n^{t}(F)$,
\bitem 
compute the group $H^{n+1}_m(F)$ in two ``classical'' cases where
$F$ is either the rational function field in one variable $F=k(t)$ 
or the formal power series $F=k((t))$.

\smallskip

We shall consider in detail the 
case  (so called ``non-classical case'')  
of  a field $F$ of characteristic $p$ and $m=p$. 


\HH 2.1. Definition of $H_{m}^{n+1}(F)$ and 
pairing $(\,\,,\,\,)_{m}$

To define the group  $H^{n+1}_m(F)$ we consider three cases 
depending on the characteristic of the field $F$.

\eg Case 1 {{\rm(Classical)}} 

 Either $\chr (F)=0$ or $\chr (F)=p$ is prime to $m$.

In this case we set 
$$
H_{m}^{n+1}(F):=H^{n+1}(F,\mu_m^{\otimes n}).
$$
The Kummer theory gives rise to the well known natural isomorphism
$F^*/F^{*m}\to H^1(F, \mu _{m})$.
Denote the image of an element $a\in F^*$ under this isomorphism  
 by $(a)$.
The cup product gives the homomorphism
$$
\underbrace{F^*\otimes\dots\otimes F^*}_{n}\to H^n(F,\mu_m^{\otimes n}),
\qquad
a_1\otimes\dots\otimes a_n\to (a_1,\dots,a_n)
$$
where $(a_1,\dots,a_n):=(a_1)\cup\dots\cup (a_n)$.
It is well known that the element $(a_1,\dots,a_n)$ is zero if $a_i+a_j=1$ for some $i\ne j$. From the definition of the Milnor $K$-group we get 
the homomorphism 
$$
\eta_m\colon
K_n^M(F)/m\to H^n(F,\mu_m^{\otimes n}),
\qquad
\{a_1,\dots, a_n\}\to (a_1,\dots,a_n).
$$
Now, we define the pairing $(\,\,,\,\,)_m$ 
as the following composite
$$
H^{1}(F, \Bbb Z/m)\otimes K_{n}(F)/m
@>{\id\otimes\eta_m}>> 
H^{1}(F, \Bbb Z/m)\otimes H^n(F,\mu_m^{\otimes n})
@>{\cup}>>
H_{m}^{n+1}(F,\mu_m^{\otimes n}).
$$
 \endeg

\eg Case 2 

$\chr(F)=p\ne 0$ 
and $m$ is a power of $p$.

To simplify the exposition we start with   
the case $m=p$.
Set 
$$
H_{p}^{n+1}(F)=\coker\bigl( \Omega _{F}^{n}
@>{\wp}>> \Omega_{F}^{n}/d\Omega _{F}^{n-1}\bigr)
$$
where  
$$
\aligned 
&d\left( adb_{2}\wedge \dots \wedge db_{n}\right)
=da\wedge db_{2}\wedge \dots \wedge db_{n}, \\
&\wp \bigl( a\frac{db_{1}}{b_{1}}\wedge 
   \dots \wedge \frac{db_{n}}{b_{n}}\bigr) 
 =\left( a^{p}-a\right) \frac{db_{1}}{b_{1}}\wedge \dots \wedge 
\frac{db_{n}}{b_{n}} + d\Omega _{F}^{n-1} 
\endaligned 
$$
($\wp=\Car^{-1}-1$ where $\Car^{-1}$ is the inverse Cartier operator
defined in subsection 4.2). 
The pairing $(\,\,,\,\,)_p$ is defined as follows:
$$
\aligned 
&(\,\,,\,\,)_{p}\colon F/\wp (F) \times  K_{n}(F)/p  
 \rightarrow  H_{p}^{n+1}(F),\\ 
&(a,\left\{b_{1},\dots ,b_{n}\right\})  
\mapsto  
a\frac{db_{1}}{b_{1}}\wedge \dots \wedge \frac{db_{n}}{b_{n}} 
\endaligned
$$
where $F/\wp (F)$ is identified with $H^{1}(F, \Bbb Z/p)$ via 
 Artin--Schreier theory.

To define the group $H^{n+1}_{p^i}(F)$ for an arbitrary $i\ge1$ 
we note that the group $H^{n+1}_p(F)$ is the 
quotient group of $\Omega_F^n$.
In particular, generators of the group $H^{n+1}_p(F)$ can be 
written in the form $ad b_1\wedge\dots\wedge db_n$.
Clearly, the natural homomorphism
$$
F\otimes\underbrace{F^*\otimes\dots\otimes F^*}_n
\to H^{n+1}_p(F),
\qquad a\otimes b_1\otimes \dots\otimes b_n\mapsto  
a\frac{db_{1}}{b_{1}}\wedge    \dots \wedge \frac{db_{n}}{b_{n}}
$$
is  surjective.
Therefore the group $H^{n+1}_p(F)$ is naturally identified with
the quotient group $F\otimes F^*\otimes\dots\otimes F^*/J$.
It is not difficult to show that the subgroup $J$
is generated by the following
elements:
\Roster 
\Item{} $(a^p-a)\otimes b_1\otimes \dots\otimes b_n$,
 
\Item{} $a\otimes a \otimes b_2\otimes\dots\otimes b_n$,

\Item{} $a\otimes b_1\otimes\dots\otimes b_n$, where $b_i=b_j$ for some
$i\ne j$.
\endRoster 
This description of the group $H^{n+1}_p(F)$ can be easily generalized to
 define   $H^{n+1}_{p^i}(F)$ for an arbitrary $i\ge1$.
Namely, we define the group 
$H^{n+1}_{p^i}(F)$ as the quotient group
$$
W_i(F)\otimes\underbrace{F^*\otimes\dots\otimes F^*}_n/J
$$
where $W_i(F)$ is the group of Witt vectors of length $i$ and
$J$ is the subgroup of 
$W_i(F)\otimes F^*\otimes\dots\otimes F^*$
generated by the following elements:
\Roster
\Item{} $({\bold F}(w)-w)\otimes b_1\otimes \dots\otimes b_n$,
\Item{} $(a,0,\dots,0)\otimes a \otimes b_2\otimes\dots\otimes b_n$, 
\Item{} $w\otimes b_1\otimes\dots\otimes b_n$, where $b_i=b_j$ for some
$i\ne j$.
\endRoster 

The pairing $(\,\,,\,\,)_{p^i}$ is defined as follows:
$$
\aligned
&(\,\,,\,\,)_{p}\colon 
W_i(F)/\wp (W_i(F)) \times  K_{n}(F)/p^i  
 \rightarrow  H_{p^i}^{n+1}(F),\\ 
&(w,\left\{b_{1},\dots ,b_{n}\right\})  
\mapsto  
w\otimes b_1\otimes\dots \otimes b_n 
\endaligned 
$$
where
$\wp={\bold F}-\id\colon W_i(F)\to W_i(F)$ and the group
$W_i(F)/\wp (W_i(F))$ is identified with $H^{1}(F, \Bbb Z/p^i)$ via 
Witt theory.
This completes definitions in Case 2.
\endeg

\eg Case 3  

$\chr (F)=p\ne 0$ and 
$m=m' p^i$ where $m'>1$ is an integer prime to $p$ and $i\ge 1$. 

\phantom{}\par 
The groups $H^{n+1}_{m'}(F)$ and $H^{n+1}_{p^i}(F)$ 
are already defined (see Cases 1 and 2).
We define the group
$H^{n+1}_m(F)$ by the following formula:
$$
H^{n+1}_m(F):=H^{n+1}_{m'}(F)\oplus H^{n+1}_{p^i}(F)
$$
Since $H^1(F,\Bbb Z/m)\simeq H^1(F,\Bbb Z/m')\oplus H^1(F,\Bbb Z/p^i)$ and
$K_n(F)/m\simeq K_n(F)/m'\oplus K_n(F)/p^i$, we can define the pairing
$(\,\,,\,\,)_m$ as the direct sum of the pairings
$(\,\,,\,\,)_{m'}$ and $(\,\,,\,\,)_{p^i}$.
This completes the definition of the group
$H^{n+1}_m(F)$ and of the pairing $(\,\,,\,\,)_m$.
 \endeg

\rk  Remark 1

In the case $n=1$ or $n=2$ the group $H^n_m(F)$ can be determined 
as follows:
$$
H^1_m(F)\simeq H^1(F,\Bbb Z/m)
\qquad{\text{and} }\qquad
H^2_m(F)\simeq {}_m\Br(F).
$$
\endrk

\rk  Remark 2

The group
$H^{n+1}_m(F)$ is often denoted by $H^{n+1}(F, \Bbb Z/m\,(n))$.
\endrk

\HH 2.2. The group $H^{n+1}(F)$

In the previous subsection we defined the group $H^{n+1}_m(F)$ 
and the pairing $(\,\,,\,\,)_m$ for an arbitrary $m$.
Now, let $m$ and $m'$ be positive integers such that $m'$ is divisible
by $m$. In this case there exists a canonical homomorphism
$$
i_{m,m'}\colon H^{n+1}_m(F)\to H^{n+1}_{m'}(F).
$$
To define the homomorphism $i_{m,m'}$ it suffices to consider the following
two cases:

\eg Case 1   

Either $\chr(F)=0$ or $\chr(F)=p$ is prime to $m$ and $m'$.

This case corresponds to Case 1 in the definition of the
group $H^{n+1}_m(F)$ (see subsection 2.1).
We identify the homomorphism $i_{m,m'}$ 
with the homomorphism 
$$
H^{n+1}(F,\mu_m^{\otimes n})\to H^{n+1}(F,\mu_{m'}^{\otimes n})
$$ 
induced by the natural embedding $\mu_m\subset \mu_{m'}$.
\endeg

\eg Case 2  

$m$ and $m'$ are powers of $p=\chr(F)$.

We can assume that $m=p^i$ and $m'=p^{i'}$ with $i\le i'$.
This case corresponds to  Case 2 in the definition of the
group $H^{n+1}_m(F)$.
We define $i_{m,m'}$ as the homomorphism induced by 
$$
\gather
W_i(F)\otimes F^*\otimes\dots F^*
\to
W_{i'}(F)\otimes F^*\otimes\dots F^*,\\
(a_1,\dots,a_i)\otimes b_1\otimes\dots\otimes b_n
\mapsto
(0,\dots,0,a_1,\dots,a_i)\otimes b_1\otimes\dots\otimes b_n.
\endgather
$$
\endeg

The maps $i_{m,m'}$ (where $m$ and $m'$ run over all integers such that
$m'$ is divisible by $m$) determine the inductive system of the groups.

\df Definition

For a field $F$ and an integer $n$ set
$$H^{n+1}(F)=\inlim_m \, H^{n+1}_m(F).$$
\enddf

\th  Conjecture 1

The natural homomorphism $H^{n+1}_m(F)\to H^{n+1}(F)$ is injective
and the image of this homomorphism coincides with the $m$-torsion part
of the group $H^{n+1}(F)$.
\endth 

This conjecture follows easily from  the Milnor--Bloch--Kato 
 conjecture (see subsection 4.1) in 
degree $n$. In particular, it is proved for  $n\le 2$.
For fields of characteristic $p$ we have the following theorem.

\th  Theorem 1 

Conjecture 1 
is true if  $\chr (F)=p$ and $m=p^i$.
\endth

\HH 2.3. Computing the group $H_{m}^{n+1}(F)$ for some fields

We start with the following well known result.
\nopagebreak

\th 
Theorem 2 {{\rm(classical)}}
\par 
Let $F$ be a perfect field.
Suppose  that $\chr(F)=0$ or $\chr(F)$ is prime to $m$. 
Then 
$$ 
\aligned
H_{m}^{n+1}\bigl( F\left( \left( t\right) \right) \bigr) \simeq
&\ H_{m}^{n+1}(F)\oplus H_{m}^{n}(F) \\
H_{m}^{n+1}\bigl( F\left( t\right) \bigr) \simeq &\ H_{m}^{n+1}(F)\oplus 
\coprod_{\text{monic irred $f(t)$}} 
H^{n}_m\left( F\left[ t\right] /f(t)\right).
\endaligned 
$$ 
\endth

It is known that we cannot omit the conditions on $F$ and $m$
in the statement of Theorem~2.
To generalize the theorem to the arbitrary case we need the following 
notation. For a complete discrete valuation field $K$
and its  maximal unramified extension   $K_{\ur}$
 define the groups $H^n_{m,\ur}(K)$ and $\~H^n_{m}(K)$ as follows:
$$
H^n_{m,\ur}(K)=\ker\,(H^n_m(K)\to H^n_m(K_{\ur}))
\qquad\hbox{ and }\quad
\~ H^n_m(K)=H^n_{m}(K)/H^n_{m,\ur}(K).
$$  
Note that  for a field $K=F((t))$ we obviously have $K_{\ur}=F^{\sep}((t))$.
We also note that under the hypotheses of Theorem~2 we have
$H^n(K)=H^n_{m,\ur}(K)$ and $H^n(K)=0$.
The following  theorem is due to Kato. 

\th Theorem 3 {{\rm(Kato, \cite{K1, Th. 3 \S 0})}}

Let $K$ be a complete discrete valuation field with 
residue field $k$.
Then 
$$
H^{n+1}_{m,\ur}(K)\simeq H^{n+1}_m(k)\oplus H^{n}_m(k).
$$
In particular, 
$H_{m,\ur}^{n+1}\left( F\left( \left( t\right) \right) \right)
\simeq H_{m}^{n+1}(F)\oplus H_{m}^{n}(F)$.
\endth

This theorem plays a key role in Kato's approach to  class field 
theory of multidimensional local fields
(see section~5 of this part). 

To generalize the second isomorphism of 
Theorem 2  we need the following  notation. Set 
$$
\aligned
&H_{m,\sep}^{n+1}(F(t))=\ker\,(H_m^{n+1}(F(t))\to H_m^{n+1}(F^{\sep}(t)))
 \text{  and  } \\
&\~ H_{m}^{n+1}(F(t))=H_{m}^{n+1}(F(t))/H_{m,\sep}^{n+1}(F(t)).
\endaligned
$$
If the field $F$ satisfies the hypotheses of Theorem 2,  
we have 

\noindent $H_{m,\sep}^{n+1}(F(t))=H_{m}^{n+1}(F(t))$
and 
$\~ H_{m}^{n+1}(F(t))=0$.

In the general case we have the following statement.

\th Theorem 4 {{\rm (Izhboldin, \cite{I2, Introduction})}}  

$$
\gather 
H_{m,\sep}^{n+1}\left( F\left( t\right) \right) \simeq
H_{m}^{n+1}(F)\oplus 
\coprod_{\text{monic irred $f(t)$}}
H^{n}_m\left( F[t]/f(t)\right), \\
\~ H_{m}^{n+1}\left( F\left( t\right) \right) 
\simeq
\coprod_{v}
\~H_{m}^{n+1}(F(t)_v)  
\endgather
$$
where $v$ runs over all normalized discrete valuations of the 
field $F(t)$ and $F(t)_v$ denotes the $v$-comple\-ti\-on 
of $F(t)$.
\endth

\HH 2.4. On the group $K_{n}(F)$

In this subsection we discuss the structure of the torsion and 
cotorsion in  Milnor $K$-theory. 
For simplicity, we consider the 
case of prime  $m=p$.
We start with the following 
fundamental theorem concerning the quotient group
$K_n(F)/p$ for fields of characteristic $p$.

\th 
Theorem 5 {{\rm (Bloch--Kato--Gabber, \cite{BK, Th. 2.1})}}

Let $F$ be a field of characteristic $p$.
Then the differential symbol 
$$ 
d_F\colon K_{n}(F)/p \to  \Omega
_{F}^{n}, \qquad 
\left\{ a_{1},\dots ,a_{n}\right\} 
\mapsto \frac{da_{1}}{a_{1}}\wedge
\dots \wedge \frac{da_{n}}{a_{n}}
 $$ 
is injective and its image  coincides with the
kernel $\nu_n(F)$ of the homomorphism
$\wp$ {{\rm(}}for the definition see Case~2 of 2.1{{\rm)}}.
In other words, the sequence
$$
\CD
0@>>> K_n(F)/p @>d_F>> \Omega_F^n @>\wp>> \Omega_F^n/d\Omega_F^{n-1}\\
\endCD
$$
is exact.
\endth

This theorem relates  the Milnor $K$-group modulo $p$ of a field of
characteristic $p$ with
 a submodule of the differential module
whose structure is  easier to understand. 
The theorem 
is important for  Kato's approach to
higher local class field theory. For a sketch of its proof see subsection~A2
in the appendix to this section.

There exists a natural generalization  of the above theorem for
the quotient groups $K_n(F)/p^i$ by using De Rham--Witt complex 
(\cite{BK, Cor. 2.8}).


\smallskip
Now, we recall well known Tate's 
conjecture concerning the torsion subgroup
of the Milnor $K$-groups.

\th Conjecture 2 {{\rm (Tate)}} 

Let $F$ be a field and $p$ be a prime integer.
\Roster 
\Item{(i)} 
If $\chr(F)\neq p$ and $\zeta_{p}\in F$, 
then ${}_{p}K_{n}(F)=\left\{ \zeta _{p}\right\} \cdot K_{n-1}(F)$. 

\Item{(ii)} If $\chr(F)=p$ then ${}_{p}K_{n}(F)=0$.
\endRoster 
\endth

\smallskip
\noindent
This conjecture is trivial in the case where $n\le1$.
In the other cases we have the following theorem.

\th  Theorem 6 

Let $F$ be a field and $n$ be a positive integer.
\Roster
\Item{(1)} 
Tate's Conjecture holds if $n\le2$  
{{\rm(}}Suslin, \cite{S}{{\rm)}},
\Item{(2)} 
Part {\rm(ii)} of Tate's Conjecture holds for all $n$ 
{\rm(}Izhboldin, \cite{I1}{\rm)}.
\endRoster 
\endth

The proof of this theorem is closely related to
the proof of Satz 90 for $K$-groups. 
Let us recall two basic conjectures on this subject.

\th Conjecture 3 {{\rm (Satz 90 for $K_n$)}}

If $L/F$ is a cyclic extension of  degree $p$ 
with the Galois group $G=\left\langle \sigma \right\rangle $ 
then the sequence
$$
K_{n}(L) @>{1-\sigma }>> K_{n}(L) @>{N_{L/F}}>> K_{n}(F)
$$
is exact.  
\endth

There is an analogue of the above conjecture for the quotient group
$K_n(F)/p$. Fix the following notation till 
the end of this section:

\df Definition 

For a field $F$ set $$k_{n}(F)=K_{n}(F)/p.$$ 
\enddf

\th Conjecture 4 {{\rm (Small Satz  90 for $k_n$)}}

If $L/F$ is a cyclic extension of  degree $p$  
with the Galois group $G=\left\langle \sigma \right\rangle $, 
then 
the sequence %
$$
k_{n}(F)\oplus k_{n}(L) @>{i_{F/L}\oplus (1-\sigma )}>>%
k_{n}(L) @>{N_{L/F}}>> k_{n}(F)
$$
is exact.
\endth 

The conjectures 2,3 and 4 are not independent:

\th  Lemma {{\rm(Suslin)}}

Fix a prime integer $p$ and integer $n$.
Then in the category of all fields {{\rm(}}of a given characteristic{{\rm)}} we have
$$
\text
{{\rm(}Small Satz 90 for $k_n${\rm)} + 
{\rm(}Tate conjecture for ${}_{p}K_n${\rm)} $\iff$ 
{\rm(}Satz 90 for $K_n${\rm)}.}
$$
Moreover, for a given field $F$ we have
$$
\text
{{\rm(}Small Satz 90 for $k_n${\rm)} + 
{\rm(}Tate conjecture for $_{p}K_n${\rm)} $\Rightarrow$
{\rm(}Satz 90 for $K_n${\rm)} }
$$
and 
$$
\text {{\rm(}Satz 90 for $K_n${\rm)} $\Rightarrow$ 
{\rm(}small Satz 90 for $k_n${\rm)}.}
$$
\endth

%

 Satz 90  conjectures are proved for  
$n\le 2$ (Merkurev-Suslin, \cite{MS1}). If $p=2$, $n=3$,
and $\chr (F)\ne 2$, the conjectures were proved by Merkurev and 
Suslin \cite{MS} and Rost. 
For $p=2$ the conjectures follow from recent results of Voevodsky. 
For fields of characteristic $p$  
the conjectures are proved for all $n$:

\th  Theorem 7 {{\rm (Izhboldin, \cite{I1})}}

Let $F$ be a field of characteristic $p$ and 
$L/F$ be a cyclic extension  of degree $p$. 
Then the following sequence is exact:
$$
0\rightarrow K_{n}(F)\rightarrow K_{n}(L) @>{1-\sigma }>>%
K_{n}(L) @>{N_{L/F}}>> K_{n}(F) \rightarrow
H^{n+1}_p(F)\rightarrow H^{n+1}_p(L)
$$ 
\endth

\HH 2.5. On the group $K_n^{t}(F)$

In this subsection we discuss the same issues, as in the previous
subsection, for the group $K_n^{t}(F)$.

\df  Definition

Let $F$ be a field and $p$ be a prime integer.
We set 
$$
DK_n(F)=\bigcap_{m\ge 1} mK_{n}(F)\quad\text{and} \quad 
D_pK_n(F)=\bigcap_{i \ge 0}p^{i}K_{n}(F).
$$
We define the group $K_n^{t}(F)$ as the quotient group:
$$
K_n^{t}(F)=K_n(F)/DK_n(F)=K_n(F)/\bigcap_{m\ge 1} mK_{n}(F).
$$
\enddf

The group $K_n^{t}(F)$ is of special interest
for higher class field theory (see sections 6, 7 and 10).
We have 
the following evident isomorphism (see also 2.0):
$$
K_n^{t}(F)\simeq \im \left(K_{n}(F) \to \prlim_{m} \, K_{n}(F)/m\right) .
$$
The quotient group $K_n^{t}(F)/m$ is obviously isomorphic to the group
$K_n(F)/m$.  
As for the torsion subgroup of $K_n^{t}(F)$, it is quite natural to state
the same questions as for the group $K_n(F)$.

\rk  Question 1

Are the $K^{t}$-analogue of Tate's conjecture and 
Satz 90 Conjecture  
true  for the group $K_{n}^{t}(F)$?
\endrk 

If we know the (positive) answer to the corresponding  
question for the group $K_n(F)$, then 
the previous question 
is  equivalent to the following:
 
\rk  Question 2

Is the group $DK_{n}(F)$ 
divisible?
\endrk

At first sight this question looks  trivial because
the group  $DK_{n}(F)$  consists of all divisible elements of $K_n(F)$.
However, the following theorem shows that the group 
$DK_n(F)$ is not necessarily a  divisible group!

\th Theorem 8 {{\rm (Izhboldin, \cite{I3})}}

For every $n\ge 2$ and prime $p$ there is a field $F$
such that  $\chr(F)\neq p$, $\zeta_{p}\in F$ and
\Roster
\Item{(1)}
The group $DK_n(F)$ is not divisible, and the group 
$D_{p}K_{2}(F)$ is not $p$-divisible,
\Item{(2)} The $K^{t}$-analogue of Tate's conjecture 
is false for  $K^{t}_n$:
$$
{}_{p} K_{n}^{t}(F)\neq \left\{ \zeta _{p}\right\} \cdot
K_{n-1}^{t}(F).
$$ 
\Item{(3)} The $K^{t}$-analogue of Hilbert 90 conjecture  
is false for  group $K^{t}_n(F)$.
\endRoster
\endth

\rk  Remark 1   

The field $F$ satisfying the conditions
of Theorem 8 can be constructed as the function field
of some  infinite dimensional variety over any field of 
characteristic zero whose group of roots of unity is finite.

Quite a different construction for irregular prime numbers $p$ 
and $F=\Bbb Q(\mu_p)$ follows from works of G. Banaszak \cite{B}. 
\endrk

\rk Remark 2

 If $F$ is a field of characteristic $p$ then the groups
$D_pK_n(F)$ and $DK_n(F)$ are $p$-divisible.
This easily implies that ${}_p K_n^{t}(F)=0$.
Moreover,  Satz 90 theorem holds for $K_n^{t}$ 
in  the case of cyclic $p$-extensions.
\endrk

\rk Remark 3

If $F$ is a multidimensional local fields then the group $K^{t}_n(F)$ 
is studied  in  section~6 of this volume. In particular, Fesenko 
(see subsections~6.3--6.8 of section~6) 
gives  
positive answers to Questions~1~and~2   
for multidimensional local fields.
\endrk

\Bib References

\rf{B}
G. Banaszak, 
{Generalization of the Moore exact sequence
and the wild kernel
for higher $K$-groups},
Compos. Math., 86(1993), 281--305.

\rf{BK}
S. Bloch and K. Kato, 
{$p$-adic \'etale cohomology},  
Inst. Hautes \'Etudes Sci. Publ. Math.  {63},
(1986), 107--152.

\rf{F}
 I. Fesenko, 
Topological Milnor $K$-groups of higher local fields,
section~6 of this volume.
 
\rf{I1}
O. Izhboldin, 
{On $p$-torsion in $K\sp M\sb *$ for fields of characteristic $p$},  
 Adv. Soviet Math., vol. 4,  129--144, Amer. Math. Soc., Providence RI, 1991

\rf{I2}
O. Izhboldin, 
On the cohomology groups of the field of rational functions,  
Mathematics in St.Petersburg, 21--44,  Amer. Math. Soc. Transl. Ser. 2, vol. {174}, 
Amer. Math. Soc., Providence, RI, 1996.

\rf{I3}
O. Izhboldin, 
On the quotient group of $K_2(F)$, 
preprint, \linebreak  
www.maths.nott.ac.uk/personal/ibf/stqk.ps 

\rf{K1}
K. Kato, 
Galois cohomology of complete discrete valuation fields, 
In {Algebraic $K$-theory}, Lect. Notes in Math. 967, Springer-Verlag, Berlin, 1982, 215--238.

\rf{K2}
K. Kato,  
{Symmetric bilinear forms, quadratic forms and Milnor $K$-theory 
    in characteristic two}, 
    Invent. Math. {66}(1982),  493--510. 

\rf{MS1}
A. S.  Merkur'ev  and A. A. Suslin, 
 {$K$-cohomology of Severi-Brauer varieties and the norm 
 residue homomorphism}, 
Izv. Akad. Nauk SSSR  
 Ser. Mat. {46}(1982);
English translation in 
 Math. USSR Izv. {21}(1983), 307--340.

\rf{MS2}
A. S.  Merkur'ev  and A. A. Suslin, 
 {The norm residue homomorphism of degree three,} 
Izv. Akad. Nauk SSSR Ser. Mat. {54}(1990);
English translation in 
 Math. USSR Izv. {36}(1991), 349--367. 

\rf{MS3}
A. S.  Merkur'ev and A. A. Suslin, 
{The group $K\sb 3$ for a field,}
Izv. Akad. Nauk SSSR Ser. Mat. { 54}(1990);
English translation in 
 Math. USSR Izv. {36}(1991), 541--565.
 
\rf{S} 
A. A.  Suslin,  {Torsion in $K_2$ of fields}, K-theory {1}(1987), 5--29. 
\endBib






\vfill
\pagebreak
\end

%% file: m3-macs.tex
\expandafter\ifx\csname mthreemacsloaded\endcsname\relax\else \fi

\magnification1100
\input amstex


 \catcode`\@=11
 \let\wlog@ld\wlog
 \def\wlog#1{\relax}

 \newif\ifIN@
 \def\m@rker{\m@@rker}
 \def\IN@{\expandafter\INN@\expandafter}
 \long\def\INN@0#1@#2@{\long\def\NI@##1#1##2##3\ENDNI@
    {\ifx\m@rker##2\IN@false\else\IN@true\fi}%
     \expandafter\NI@#2@@#1\m@rker\ENDNI@}
  \newtoks\Initialtoks@  \newtoks\Terminaltoks@
  \def\SPLIT@{\expandafter\SPLITT@\expandafter}
  \def\SPLITT@0#1@#2@{\def\TTILPS@##1#1##2@{%
     \Initialtoks@{##1}\Terminaltoks@{##2}}\expandafter\TTILPS@#2@}
  \newtoks\Trimtoks@

 \def\ForeTrim@{\expandafter\ForeTrim@@\expandafter}
 \def\ForePrim@0 #1@{\Trimtoks@{#1}}
 \def\ForeTrim@@0#1@{\IN@0\m@rker. @\m@rker.#1@%
     \ifIN@\ForePrim@0#1@%
     \else\Trimtoks@\expandafter{#1}\fi}
 
  \def\Trim@0#1@{%
      \ForeTrim@0#1@%
      \IN@0 @\the\Trimtoks@ @%
        \ifIN@
             \SPLIT@0 @\the\Trimtoks@ @\Trimtoks@\Initialtoks@
             \IN@0\the\Terminaltoks@ @ @%
                 \ifIN@
                 \else \Trimtoks@ {FigNameWithSpace}%
                 \fi
        \fi
      }

  \font\titlebold=cmbx12 scaled 1200
  \font\twelvebold=cmbx12
  \font\tenbold=cmbx10
  \font\ninebold=cmbx9
  \font\sevenbold=cmbx7
  \font\fivebold=cmbx5

  \input amssym.def \input amssym
     \font\titlemsa=msam10 at 14.4pt
     \font\titlemsb=msbm10 at 14.4pt
     \font\titleeufm=eufm10 at 14.4pt
     \font\twelvemsa=msam10 scaled 1200
     \font\twelvemsb=msbm10 scaled 1200
     \font\twelveeufm=eufm10 scaled 1200
     \font\ninemsa=msam9
     \font\ninemsb=msbm9
     \font\nineeufm=eufm9

   \ifx\cyrfam\undefined
   \else
     \immediate\write16{}%
     \message{ !!! cyr fonts already defined. !!! }
     \message{ --- edit out superfluous font defs? }
   \fi
   \newfam\cyrfam
       \font\titlecyr=wncyr10 scaled 1440 
       \font\twelvecyr=wncyr10 scaled 1200
       \font\tencyr=wncyr10
       \font\ninecyr=wncyr9
       \font\sevencyr=wncyr7
       \font\sixcyr=wncyr6

   \newfam\eusmfam
       \font\titleeusm=eusm10 scaled 1440
       \font\twelveeusm=eusm10 scaled 1200
       \font\teneusm=eusm10
       \font\nineeusm=eusm9
       \font\seveneusm=eusm7
       
       \font\fiveeusm=eusm5

    \font\ninemrm=cmr9 
    \font\ninei=cmmi9
    \font\ninesy=cmsy9 
    \skewchar\ninei='177
    \skewchar\ninesy='60

  \font\twelvemrm=cmr10 at 12pt 
  \font\twelvei=cmmi10 at 12pt
  \font\twelvesy=cmsy10 at 12pt

  \font\titlemrm=cmr10 at 14.4pt 
  \font\titlei=cmmi10 at 14.4pt
  \font\titlesy=cmsy10 at 14.4pt


  \def\Smallfonts{\ninepoint}

  \def\Hfont{\titlepoint\bf}
  \def\Authorfont{\twelvepoint\it}
  \def\HHfont{\twelvepoint\bf}
  \def\HHHfont{\bf}
  \def\Bibfont{\tenbf}
  \def\Coordfont{\nineit }

  \def \thfont {\bf }
  \def \pffont {\it\itSpacing }
  \def \rkfont {\bf }
  \def \dffont {\bf }
  \def \egfont {\bf }

 \def\ninepoint{%
  \def\rm{\fam0\ninerm}%
    \textfont0=\ninemrm  \scriptfont0=\sevenrm  \scriptscriptfont0=\fiverm
    \textfont1=\ninei    \scriptfont1=\seveni   \scriptscriptfont1=\fivei
  \def\mit{\fam1\ninei}%
  \def\oldstyle{\fam1\ninei}%
    \textfont2=\ninesy   \scriptfont2=\sevensy  \scriptscriptfont2=\fivesy
    \textfont3=\tenex    \scriptfont3=\tenex    \scriptscriptfont3=\tenex
  \def\it{\fam\itfam\nineit}%
    \textfont\itfam=\nineit
  \def\bf{\ifmmode\fam\bffam\else\ninebf\fi}%
    \textfont\bffam=\ninebold 
    \scriptfont\bffam=\sevenbold 
    \scriptscriptfont\bffam=\fivebold%
  \def\msa{\fam\msafam\ninemsa}%
    \textfont\msafam=\ninemsa 
    \scriptfont\msafam=\sevenmsa
    \scriptscriptfont\msafam=\fivemsa%
  \def\msb{\fam\msbfam\ninemsb}%
    \textfont\msbfam=\ninemsb%
    \scriptfont\msbfam=\sevenmsb%
    \scriptscriptfont\msbfam=\fivemsb%
  \def\eufm{\fam\eufmfam\nineeufm}%
    \textfont\eufmfam=\nineeufm
    \scriptfont\eufmfam=\seveneufm
    \scriptscriptfont\eufmfam=\fiveeufm
   \def\eusm{\fam\eusmfam\nineeusm}%
     \textfont\eusmfam=\nineeusm
     \scriptfont\eusmfam=\seveneusm
     \scriptscriptfont\eusmfam=\fiveeusm
   \def\cyr{\fam\cyrfam\ninecyr}%
     \textfont\cyrfam=\ninecyr
     \scriptfont\cyrfam=\sevencyr
     \scriptscriptfont\cyrfam=\sixcyr
  \setbox\strutbox=\hbox{\vrule
      height7pt depth3pt width0pt}%
   \baselineskip=10.8pt\rm}

 \let\eightpoint\ninepoint 

 \def\tenpoint{%
  \def\rm{\fam0\tenrm}%
    \textfont0=\tenmrm \scriptfont0=\sevenrm \scriptscriptfont0=\fiverm%
  \def\mit{\fam1\teni}%
  \def\oldstyle{\fam1\teni}%
    \textfont1=\teni   \scriptfont1=\seveni  \scriptscriptfont1=\fivei%
    \textfont2=\tensy  \scriptfont2=\sevensy \scriptscriptfont2=\fivesy%
    \textfont3=\tenex  \scriptfont3=\tenex   \scriptscriptfont3=\tenex%
  \def\it{\fam\itfam\tenit}%
    \textfont\itfam=\tenit%
  \def\bf{\ifmmode\fam\bffam\else\tenbf\fi}%
    \textfont\bffam=\tenbold
    \scriptfont\bffam=\sevenbold%
    \scriptscriptfont\bffam=\fivebold%
  \def\msa{\fam\msafam\tenmsa}%
    \textfont\msafam=\tenmsa%
    \scriptfont\msafam=\sevenmsa%
    \scriptscriptfont\msafam=\fivemsa%
  \def\msb{\fam\msbfam\tenmsb}%
    \textfont\msbfam=\tenmsb%
    \scriptfont\msbfam=\sevenmsb%
    \scriptscriptfont\msbfam=\fivemsb%
  \def\eufm{\fam\eufmfam\teneufm}%
   \textfont\eufmfam=\teneufm
   \scriptfont\eufmfam=\seveneufm
   \scriptscriptfont\eufmfam=\fiveeufm
   \def\eusm{\fam\eusmfam\teneusm}%
    \textfont\eusmfam=\teneusm
    \scriptfont\eusmfam=\seveneusm
    \scriptscriptfont\eusmfam=\fiveeusm
   \def\cyr{\fam\cyrfam\tencyr}%
    \textfont\cyrfam=\tencyr
    \scriptfont\cyrfam=\sevencyr
    \scriptscriptfont\cyrfam=\sixcyr
  \setbox\strutbox=\hbox{\vrule %
      height8.5pt depth3.5ptwidth0pt}%
  \baselineskip=\StdBaselineskip\rm}

 \def\twelvepoint{%
  \def\rm{\fam0\twelverm}%
    \textfont0=\twelvemrm \scriptfont0=\tenmrm \scriptscriptfont0=\sevenrm
    \textfont1=\twelvei   \scriptfont1=\teni   \scriptscriptfont1=\seveni
  \def\mit{\fam1\twelvei}%
  \def\oldstyle{\fam1\twelvei}%
    \textfont2=\twelvesy  \scriptfont2=\tensy  \scriptscriptfont2=\sevensy
    \textfont3=\tenex  \scriptfont3=\tenex  \scriptscriptfont3=\tenex
  \def\it{\fam\itfam\twelveit}%
    \textfont\itfam=\twelveit
  \def\bf{\ifmmode\fam\bffam\else\twelvebf\fi}%
    \textfont\bffam=\twelvebold
    \scriptfont\bffam=\tenbold%
    \scriptscriptfont\bffam=\sevenbold%
  \def\msa{\fam\msafam\twelvemsa}%
    \textfont\msafam=\twelvemsa%
    \scriptfont\msafam=\tenmsa%
    \scriptscriptfont\msafam=\sevenmsa%
  \def\msb{\fam\msbfam\twelvemsb}%
    \textfont\msbfam=\twelvemsb%
    \scriptfont\msbfam=\tenmsb%
    \scriptscriptfont\msbfam=\sevenmsb%
  \def\eufm{\fam\eufmfam\twelveeufm}%
   \textfont\eufmfam=\twelveeufm
   \scriptfont\eufmfam=\teneufm
   \scriptscriptfont\eufmfam=\seveneufm
   \def\eusm{\fam\eusmfam\twelveeusm}%
    \textfont\eusmfam=\twelveeusm
    \scriptfont\eusmfam=\teneusm
    \scriptscriptfont\eusmfam=\seveneusm
   \def\cyr{\fam\cyrfam\tencyr}%
    \textfont\cyrfam=\twelvecyr
    \scriptfont\cyrfam=\tencyr
    \scriptscriptfont\cyrfam=\sevencyr
  \setbox\strutbox=\hbox{\vrule
      height10.2pt depth4.55pt width0pt}%
  \baselineskip=14pt\rm}

 \def\titlepoint{%
    \textfont0=\titlemrm \scriptfont0=\twelvemrm \scriptscriptfont0=\tenmrm
    \textfont1=\titlei   \scriptfont1=\twelvei   \scriptscriptfont1=\teni
  \def\mit{\fam1\titlei}%
  \def\oldstyle{\fam1\titlei}%
    \textfont2=\titlesy  \scriptfont2=\twelvesy  \scriptscriptfont2=\tensy
    \textfont3=\tenex
    \scriptfont3=\tenex
    \scriptscriptfont3=\tenex
  \def\it{\fam\itfam\titleit}%
    \textfont\itfam=\titleit
  \def\bf{\ifmmode\fam\bffam\else\titlebf\fi}%
    \textfont\bffam=\titlebold
    \scriptfont\bffam=\twelvebold%
    \scriptscriptfont\bffam=\tenbold%
  \def\msa{\fam\msafam\titlemsa}%
    \textfont\msafam=\titlemsa%
    \scriptfont\msafam=\twelvemsa%
    \scriptscriptfont\msafam=\tenmsa%
  \def\msb{\fam\msbfam\titlemsb}%
    \textfont\msbfam=\titlemsb%
    \scriptfont\msbfam=\twelvemsb%
    \scriptscriptfont\msbfam=\tenmsb%
  \def\eufm{\fam\eufmfam\titleeufm}%
    \textfont\eufmfam=\titleeufm
    \scriptfont\eufmfam=\twelveeufm
    \scriptscriptfont\eufmfam=\teneufm
   \def\eusm{\fam\eusmfam\titleeusm}%
     \textfont\eusmfam=\titleeusm
     \scriptfont\eusmfam=\twelveeusm
     \scriptscriptfont\eusmfam=\teneusm
   \def\cyr{\fam\cyrfam\tencyr}%
    \textfont\cyrfam=\titlecyr
    \scriptfont\cyrfam=\twelvecyr
    \scriptscriptfont\cyrfam=\tencyr
  \setbox\strutbox=\hbox{\vrule
      height12.3pt depth5.54pt width0pt}%
  \baselineskip=16pt\rm}

\newbox\AuthorBox\newbox\TitleBox
\newbox\TFLinebox
\newbox\FLinebox
\newbox\HLinebox
\def\SetTFLinebox#1{\setbox\TFLinebox=\hbox{#1}}
\def\SetFLinebox#1{\setbox\FLinebox=\hbox{#1}}
\def\SetHLinebox#1{\setbox\HLinebox=\hbox{#1}}

 \def\SetAuthorHead#1{%
     \setbox\AuthorBox=\hbox{\ninepoint \it 
           \ignorespaces\frenchspacing#1\unskip}}
 \def\SetTitleHead#1{%
     \setbox\TitleBox=\hbox{\ninepoint \it
           \ignorespaces\frenchspacing#1\unskip}}

  \def\itSpacing{\relax}
  \def\itSpacingOff{\relax}


 \def\Hrule{\hrule width0pt height0pt}

  \newskip\ProcSkip \ProcSkip 8pt plus2pt minus2pt

 \newskip\LastSkip
 \def\SaveLastSkip{\LastSkip\lastskip}
 \def\RestoreLastSkip{\vskip-\LastSkip\vskip\LastSkip}

 \def\NoindentAfter{\everypar={\setbox0=\lastbox\everypar={}}}

 \long\def\H#1\par#2\par{\notenumber=0 \titlepagetrue%
    {
    \baselineskip=20pt
    \parindent=0pt\parskip=0pt\frenchspacing
    \leftskip=0pt plus .2\hsize minus .3\hsize
    \rightskip=0pt plus .2\hsize minus .3\hsize
 \def\\{\unskip\break}%
    \pretolerance=10000 \Hfont #1\unskip\break
     \vskip7pt\Hrule
\hfill \Authorfont #2\hfill\hfill\unskip}
    \vskip48pt plus 4pt minus 4pt
    \par\NoindentAfter\rm}

 \long\def\Hi#1\par#2\par{\notenumber=0 \titlepagetrue%
    {  \baselineskip=0pt  \parindent=0pt\parskip=0pt\frenchspacing
    \leftskip=0pt plus .2\hsize minus .3\hsize
    \rightskip=0pt plus .2\hsize minus .3\hsize
}
    \rm}


 \newdimen\PageRemainder
  \def\SetPageRemainder{
     \PageRemainder=\pagegoal
     \ifdim\PageRemainder=\maxdimen\PageRemainder=\vsize
     \else\advance\PageRemainder by -1\pagetotal\fi}

  \def\Rpt@{}\def\Rpt@@{}

  \long\def\HH#1\par{\par
  \SaveLastSkip\removelastskip\goodbreak
  \ifdim\LastSkip<30pt 
     \LastSkip 30pt
plus 3pt minus 2pt\fi
  \SetPageRemainder\advance\PageRemainder-\LastSkip
  \ifdim\PageRemainder<150pt
       \edef\Rpt@{remain = \the\PageRemainder\noexpand\\
                pagetotal=\the\pagetotal\noexpand\\
                           pagegoal=\the\pagegoal}%
          \fi
   \ifdim\PageRemainder<65pt 
       \ifdim\PageRemainder > 0pt
          \edef\Rpt@@{\noexpand\\
                      Had HH PageRemainder$<$\relax 65pt\noexpand\\
                      Hence forced break!}%
     \vskip 0pt plus .2\PageRemainder\eject 
    \fi\fi
    \vskip\LastSkip\Hrule 
    \pretolerance=10000\rightskip=0pt plus 3em
    \hangafter1 \hangindent=2.2em%
    \noindent
    \HHfont \unskip \Ednote{\Rpt@\Rpt@@}%
            \def\Rpt@{}\def\Rpt@@{}%
            \ignorespaces
            #1\par\rightskip=0pt\pretolerance=\StdPretolerance%
    \NoindentAfter
\tenpoint\rm%
     \medskip \vskip\ProcSkip}

  \long\def\HHH#1\par{\par%
  \SaveLastSkip\removelastskip\goodbreak
  \ifdim\LastSkip<\ProcSkip%
     \LastSkip\ProcSkip\fi
  \SetPageRemainder\advance\PageRemainder-\LastSkip
  \ifdim\PageRemainder<150pt
       \edef\Rpt@{remain = \the\PageRemainder\noexpand\\
                pagetotal=\the\pagetotal\noexpand\\
                           pagegoal=\the\pagegoal}%
       \fi
   \ifdim\PageRemainder<48pt  
        \ifdim\PageRemainder > 0pt
             \edef\Rpt@@{\noexpand\\
                      Had HHH PageRemainder$<$\relax48pt\noexpand\\
                      Hence forced break!}%
       \vskip 0pt plus .2\PageRemainder\eject 
      \fi\fi
   \vskip\LastSkip\par\noindent
   \HHHfont \unskip\Ednote{\Rpt@\Rpt@@}%
  \def\Rpt@{}\def\Rpt@@{}%
  \ignorespaces
   #1\unskip.\quad\rm\ignorespaces
   \ignorepars}

  \long\def\ignorepars#1\par{\def\Test{#1}%
     \ifx\Test\Empty\def\This{\ignorepars}%
        \else\def\This{\Test\par}\fi
           \This}
  \def\Empty{}

 \def\Abstract#1\par{\bgroup\Smallfonts\narrower\HHH #1\par}
 \def\endAbstract{\par\egroup}


 \def\ProcBreak{\par%
    \ifdim\lastskip<8pt%
    \removelastskip%
    \penalty-200\vskip\ProcSkip\fi}

 \def\th#1\par{\ProcBreak \noindent
   {\thfont\ignorespaces
    #1\unskip.}\it\itSpacing\kern.4em\ignorepars}

 \def\endth{\ProcBreak\rm\itSpacingOff }


 \def\pf#1\par{\ProcBreak %
    \noindent\pffont#1\unskip.\rm\itSpacingOff{\kern .7em}\ignorepars}


  \def\qedbox{\hbox{\vbox{
    \hrule width0.2cm height0.2pt
    \hbox to 0.2cm{\vrule height 0.2cm width 0.2pt
             \hfil\vrule height0.2cm width 0.2pt}
    \hrule width0.2cm height 0.2pt}\kern1pt}}

  \def\qed{\ifmmode\qedbox
    \else\unskip\ \hglue0mm\hfill\qedbox\ProcBreak\fi}

  \def \rk #1\par{\ProcBreak
     \noindent{\rkfont\ignorespaces #1\unskip.}%
     \rm\kern.6em\ignorepars}

  \def \endrk {\medskip\ProcBreak }

  \def \df #1\par{\ProcBreak
     \noindent{\dffont\unskip\ignorespaces #1\unskip.}%
     \rm\kern.6em\ignorepars}

  \def \enddf {\medskip\ProcBreak }

  \def \eg #1\par{\ProcBreak
     \noindent\egfont\unskip\ignorespaces #1\unskip.
     \rm\kern.6em\ignorepars}

  \def \endeg {\medskip\ProcBreak }

  \newdimen\Overhang

   \def\MaxTag@#1#2#3#4#5{\setbox0=\hbox{#4\ignorespaces#2\unskip}%
     \dimen0=\wd0\advance\dimen0 by#3
     \ifdim\dimen0<#5\relax\dimen0=#5\fi
     \expandafter\edef\csname #1Hang\endcsname{\the\dimen0}}

 \def\MaxItemTag#1{\MaxTag@{Item}{#1}{.4em}{\ItemStyle}{\parindent}}%
 \def\MaxItemItemTag#1{%
        \MaxTag@{ItemItem}{#1}{.4em}{\ItemItemStyle}{\parindent}}
 \def\MaxNrTag#1{\MaxTag@{Nr}{#1}{.5em}{\NrStyle}{\parindent}}
 \def\MaxReferenceTag#1{%
        \MaxTag@{Reference}{[#1]}{.6em}{\ninerm}{\parindent}}
 \def\MaxFootTag#1{\MaxTag@{Foot}{#1}{.4em}{\ninerm}{\z@}}

  \def\SetOverhang@{\Overhang=.8\dimen0%
     \advance\Overhang by \wd0\relax
     \ifdim\Overhang>\hangindent\relax
       \advance\Overhang by .25\dimen0%
       \Ednote{Tag is pushing text.}\osumess{Tag is pushing text.}%
     \else\Overhang=\hangindent
     \fi}

   \def\Item#1{\par\noindent
      \hangafter1\hangindent=\ItemHang
      \setbox0=\hbox{\ItemStyle\ignorespaces#1\unskip}%
      \dimen0=.4em\SetOverhang@
      \rlap{\box0}\kern\Overhang\ignorespaces}

   \def\ItemItem#1{\par\noindent
      \hangafter1\hangindent=\ItemItemHang
      \setbox0=\hbox{\ItemItemStyle\ignorespaces#1\unskip}%
      \dimen0=.4em\SetOverhang@
      \advance\hangindent by \ItemHang
      \kern\ItemHang\rlap{\box0}%
      \kern\Overhang\ignorespaces}

  \def\Nr#1{\par\noindent\hangindent=\NrHang 
    \setbox0=\hbox{\NrStyle\ignorespaces#1\unskip}%
    \dimen0=.5em\SetOverhang@
    \rlap{\box0}\kern\Overhang
    \hangindent=\z@\ignorespaces}

   \newskip\Rosterskip\Rosterskip 1pt plus1pt 
   \def\Roster{\par\ifdim\lastskip<\Rosterskip\removelastskip\vskip\Rosterskip\fi
    \bgroup}
   \def\endRoster{\par\global\edef\LastSkip@{\the\lastskip}\removelastskip
       \egroup\penalty-50\LastSkip\LastSkip@\relax
       \ifdim\LastSkip<\Rosterskip\LastSkip\Rosterskip\fi
       \vskip\LastSkip}




 \def\cite#1{
    \def\nextiii@##1,##2\end@{{\frenchspacing\rm 
      \lBr\ignorespaces##1\unskip{\rm,~\ignorespaces##2}\rBr}}%
    \IN@0,@#1@%
    \ifIN@\def\next{\nextiii@#1\end@}\else
    \def\next{{\rm\lBr#1\rBr}}\fi\next}


   \def \Bib#1\par{%
       \par\removelastskip\SetPageRemainder
       \ifdim\PageRemainder < 97pt
        \ifdim\PageRemainder > 0pt
        \vfill\eject
       \fi\fi
    \ProcBreak \par\begingroup\parskip=0 pt%
    \goodbreak \vskip 15 pt plus 10 pt
    \noindent\null\hfill\Bibfont
      \ignorespaces #1\unskip\hfill\null\par 
    \frenchspacing \Smallfonts\rm
    \parskip=2.5 pt plus 1 pt minus.5pt%
    \nobreak\vskip 12pt plus 2pt minus2pt\nobreak
    \leftskip=0 pt \baselineskip=10.5pt}

 \def\ReferenceTagSlide{0em}
  \def\ReferenceTagGap{.5em}

  \def \rf#1{\par\noindent
     \hangafter1\hangindent=\ReferenceHang      
     \setbox0=\hbox{\ninerm[\ignorespaces#1\unskip]}%
     \dimen0=\ReferenceTagGap\SetOverhang@
     \rlap{\kern\ReferenceTagSlide\box0}%
     \kern\Overhang\ignorespaces}

  \def\ref#1\par#2\par#3\par#4\par{%
     \rf{#1}#2\unskip,\ #3\unskip,\
     #4\unskip.}

  \def\endBib{\par\endgroup\vskip 12pt minus 6pt }


  \long\def\Coordinates#1\endCoordinates{
 {\par\vskip4pt\def\\{\unskip, }\Coordfont\baselineskip10.5pt\noindent#1}}

 \def\pagecontents{
  \gdef\Pagetot@l{\pagetotal}
  \ifvoid\TRMargIns\else
    \rlap{\kern\hsize\kern10pt\vbox to 0pt{%
         \box\TRMargIns\vss}}\fi
  \ifvoid\topins\else\unvbox\topins\fi
   \dimen@=\dp\@cclv \unvbox\@cclv 
   \ifvoid\footins\else 
     \vskip\skip\footins
     \footnoterule
     \unvbox\footins\fi
   \ifr@ggedbottom \kern-\dimen@ \vfil \fi}


 \newcount\Ht 

 \def \Acc{\expandafter } 

 \def\swthat{\raise -1.1 ex\hbox{\sam$\widehat{}$}}
 \def\swttilde{\raise -1.2 ex\hbox{\sam$\widetilde{}$}}
 \def \overdot{{\raise .2 ex \hbox to 0pt {\hss\bf\smash{.}\hss}}}
 \def \overcircle{{\raise .1 ex \hbox to 0pt
    {\sam$\eightpoint\scriptstyle\hss\circ\hss$}}}

 \def \Mathaccent#1#2{{\sam 
  \setbox4=\hbox{$\vphantom{#2}$}
  \Ht=\ht4 
  \setbox5=\hbox{${#1}$}
  \setbox6=\hbox{${#2}$}
  \setbox7=\hbox to .5\wd6{}
  \copy7\kern .1\Ht \raise\Ht sp\hbox{\copy5}\kern-.1\Ht
  \copy7\llap{\box6}
  }}

  \def\SwtCheck #1{
        \ifmmode \check{#1}%
                \else \v {#1}%
                \fi}

 \def\barpartial {%
   \kern .17 em
    \overline {\kern -.17 em\partial\kern-.03 em}%
    \kern .03 em}

 
  \def\Overline#1{\setbox1=\hbox{\sam ${#1}$}%
      \ifdim \wd1 > 6pt
    \kern .11 em
    \overline {\kern -.11 em#1\kern-.14 em}
    \kern .14 em
  \else
    \kern .03 em
    \overline {\kern -.03 em#1\kern-.04 em}
    \kern .04 em
  \fi}

 \def\SOverline#1{\setbox1=\hbox{\sam ${#1}$}%
      \ifdim \wd1 > 7pt
    \kern .22 em
    \overline {\kern -.22 em#1\kern-.09 em}%
    \kern .09 em
  \else
    \kern .10 em
    \overline {\kern -.10 em#1\kern-.04 em}%
    \kern .04 em
  \fi}


 \def\Underline#1{\setbox1=\hbox{\sam ${#1}$}%
      \ifdim \wd1 > 6pt
    \kern .11 em
    \underline {\kern -.11 em#1\kern-.14 em}
    \kern .14 em
  \else
    \kern .03 em
    \underline {\kern -.03 em#1\kern-.04 em}
    \kern .04 em
  \fi}

 \def\SUnderline#1{\setbox1=\hbox{\sam ${#1}$}%
      \ifdim \wd1 > 7pt
    \kern .04 em
    \underline {\kern -.04 em#1\kern-.2 em}%
    \kern .2 em
  \else
    \kern .0 em
    \underline {\kern -.0 em#1\kern-.15 em}%
    \kern .15 em
  \fi}


 \def \Blackbox
   {\leavevmode\hskip .3pt \vbox
   {\hrule height 5pt\hbox{\hskip 4.5pt}}\hskip .5pt}

 \def \XX{\Blackbox\kern.5pt\Blackbox} 

  \def\.{.\kern1pt}

    \def\Hyphen{\edef\this{\the\hyphenchar\font}%
          \hyphenchar\font=-1\char\this\hyphenchar\font=\this}

 \ifx\undefined\text
  \def\text#1{\hbox{\rm #1}}\fi 



   \everymath{}  

  \def\PassMath@@{\aftergroup\AfterMath@} 

 \let\PassMath@\PassMath@@

 \def\AfterMath@{\futurelet\next\AfterMathMole@}

 \def\AfterMathMole@{
      \ifcat\next\space
          \def\this{}
      \else
      \ifcat\next\egroup %
        \def\this{\osumess{Handset mathsurround?? ---(see dollar brace)}}%
      \else
      \def\this{\AAfterMath@}
      \fi\fi
      \this}

 \def\hyphen@{-}
 \def\paren@{)}
 \def\apostr@{'}

 \def\MSC#1{\kern-.8\mathsurround#1\kern.8\mathsurround}

 \def\AAfterMath@#1{\def\Next{#1}
    \IN@0\Next @,.;:!?\relax @%
    \ifIN@\def\this{\MSC{\Next}}%
    \else
    \ifx\Next\hyphen@\def\this{\futurelet\next\AfterHyphen@}%
    \else
    \ifx\Next\paren@\def\this{#1}%
    \else 
    \ifx\Next\apostr@\def\this{#1}%
    \else \def\this{\osumess{Handset mathsurround??}%
                 #1}\fi\fi\fi\fi
    \this}

 \def\AfterHyphen@#1{\def\Next{#1}%
   \ifx\Next\hyphen@\def\this{--}\else
   \ifcat\next\space%
   \def\this{\kern-\mathsurround\kern.05em- \Next}\else
   \def\this{\kern-\mathsurround\kern.05em\Hyphen\Next}\fi\fi\this}

 \def\sam{\mathsurround=\z@\let\PassMath@\relax}  %
 \def\mas{\mathsurround=\StdMathsurround\let\PassMath@\PassMath@@}
 
 \def\Mas{\mathsurround=\StdMathsurround
                \everymath{\PassMath@}\let\PassMath@\PassMath@@}

 \def\m@th{\mathsurround=\z@\everymath{}}

 \def\m@@th{\mathsurround=\z@\everymath={}\let\m@th\relax}

\def\underbar#1{$\setbox\z@\hbox{#1}\dp\z@\z@
      \m@th \underline{\box\z@}$\relax}

\def\mathhexbox#1#2#3{\leavevmode
  \hbox{\m@@th$\m@th \mathchar"#1#2#3$}}

\def\dots{\relax\ifmmode\ldots\else$\m@th\ldots\,$\relax\fi}

\def\dotfill{\cleaders\hbox{\m@@th$\m@th \mkern1.5mu.\mkern1.5mu$}\hfill}
\def\rightarrowfill{$\m@th\mathord-\mkern-6mu%
  \cleaders\hbox{\m@@th$\mkern-2mu\mathord-\mkern-2mu$}\hfill
  \mkern-6mu\mathord\rightarrow$\relax}
\def\leftarrowfill{$\m@th\mathord\leftarrow\mkern-6mu%
  \cleaders\hbox{\m@@th$\mkern-2mu\mathord-\mkern-2mu$}\hfill
  \mkern-6mu\mathord-$\relax}

\def\downbracefill{$\m@th\braceld\leaders\vrule\hfill\braceru
  \bracelu\leaders\vrule\hfill\bracerd$\relax}
\def\upbracefill{$\m@th\bracelu\leaders\vrule\hfill\bracerd
  \braceld\leaders\vrule\hfill\braceru$\relax}

\def\angle{{\vbox{\m@@th\ialign{$\m@th\scriptstyle##$\crcr
      \not\mathrel{\mkern14mu}\crcr
      \noalign{\nointerlineskip}
      \mkern2.5mu\leaders\hrule height.34pt\hfill\mkern2.5mu\crcr}}}}

\def\big#1{{\m@@th\hbox{$\left#1\vbox to8.5\p@{}\right.\n@space$}}}
\def\Big#1{{\m@@th\hbox{$\left#1\vbox to11.5\p@{}\right.\n@space$}}}
\def\bigg#1{{\m@@th\hbox{$\left#1\vbox to14.5\p@{}\right.\n@space$}}}
\def\Bigg#1{{\m@@th\hbox{$\left#1\vbox to17.5\p@{}\right.\n@space$}}}
\def\n@space{\nulldelimiterspace\z@ \m@th}

\def\root#1\of{\setbox\rootbox\hbox{\m@@th$\m@th\scriptscriptstyle{#1}$}
  \mathpalette\r@@t}
\def\r@@t#1#2{\setbox\z@\hbox{\m@@th$\m@th#1\sqrt{#2}$\relax}
  \dimen@\ht\z@ \advance\dimen@-\dp\z@
  \mkern5mu\raise.6\dimen@\copy\rootbox \mkern-10mu \box\z@}

\def\mathph@nt#1#2{\setbox\z@\hbox{\m@@th$\m@th#1{#2}$}\finph@nt}

\def\mathsm@sh#1#2{\setbox\z@\hbox{\m@@th$\m@th#1{#2}$}\finsm@sh}

\def\@vereq#1#2{\lower.5\p@\vbox{\m@@th\baselineskip\z@skip\lineskip-.5\p@
    \ialign{$\m@th#1\hfil##\hfil$\crcr#2\crcr=\crcr}}}

\def\mathpalette#1#2{\sam\mathchoice{#1\displaystyle{#2}}%
  {#1\textstyle{#2}}{#1\scriptstyle{#2}}{#1\scriptscriptstyle{#2}}\mas}

\def\widehat#1{\setbox\z@\hbox{\sam$#1$}%
 \ifdim\wd\z@>\tw@ em\mathaccent"0\msbfam@5B{#1}%
 \else\mathaccent"0362{#1}\fi}
\def\widetilde#1{\setbox\z@\hbox{\sam$#1$}%
 \ifdim\wd\z@>\tw@ em\mathaccent"0\msbfam@5D{#1}%
 \else\mathaccent"0365{#1}\fi}

 \def\dots{\relax{}
  \ifmmode\def\thedots{\mdots@}\else\def\thedots{\tdots@}\fi %
  \thedots}

 \let\@ldeqno\eqno\let\@ldleqno\leqno
 \def\eqno{\everymath{}\@ldeqno} \def\leqno{\everymath{}\@ldleqno}

  \let\@ldeqalignno\eqalignno
  \def\eqalignno#1{\sam\@ldeqalignno{#1}\mas}
  \let\@ldeqalign\eqalign
  \def\eqalign#1{\sam\@ldeqalign{#1}\mas}

 \def\overrightarrow#1{\vbox{\m@th\ialign{##\crcr
      \rightarrowfill\crcr\noalign{\kern-\p@\nointerlineskip}
      $\hfil\displaystyle{#1}\hfil$\crcr}}}
 \def\overleftarrow#1{\vbox{\m@th\ialign{##\crcr
      \leftarrowfill\crcr\noalign{\kern-\p@\nointerlineskip}
      $\hfil\displaystyle{#1}\hfil$\crcr}}}
 \def\overbrace#1{\mathop{\vbox{\m@th\ialign{##\crcr\noalign{\kern3\p@}
      \downbracefill\crcr\noalign{\kern3\p@\nointerlineskip}
      $\hfil\displaystyle{#1}\hfil$\crcr}}}\limits}
 \def\underbrace#1{\mathop{\vtop{\m@th\ialign{##\crcr
      $\hfil\displaystyle{#1}\hfil$\crcr\noalign{\kern3\p@\nointerlineskip}
      \upbracefill\crcr\noalign{\kern3\p@}}}}\limits}

  \let\@ldmatrix\matrix
  \let\end@ldmatrix\endmatrix
  \def\matrix{\sam\@ldmatrix}
  \def\endmatrix{\end@ldmatrix\mas}
  \let\@ldgather\gather
  \let\end@ldgather\endgather
  \def\gather{\sam\@ldgather}
  \def\endgather{\end@ldgather\mas}
  \let\@ldalign\align
  \let\end@ldalign\endalign
  \def\align{\sam\@ldalign}
  \def\endalign{\end@ldalign\mas}
  \let\@ldaligned\aligned
  \let\end@ldaligned\endaligned
  \def\aligned{\sam\@ldaligned}
  \def\endaligned{\end@ldaligned\mas}
  \let\@ldtag\tag
  \def\tag{\sam\@ldtag}
   %

   \let\MinCDArrowWidth\minCDaw@




\newskip\insertskipamount\newskip\inserthardskipamount
\insertskipamount 6pt plus2pt 
\inserthardskipamount 6pt
\def\insertskip{\vskip\insertskipamount}
\newcount\SplitTest
\def\SetSplitTest{\SplitTest\insertpenalties
  \insert\topins{\floatingpenalty1}%
  \advance\SplitTest-\insertpenalties}
\def\midinsert{\par
 \SaveLastSkip\penalty-150\SetSplitTest\RestoreLastSkip
 \ifnum\SplitTest=-1
  \@midfalse\p@gefalse\else\@midtrue\fi\@ins}
\def\@ins{\par\begingroup\setbox\z@\vbox\bgroup%
  \vglue\inserthardskipamount}
\def\endinsert{\egroup 
  \if@mid \dimen@\ht\z@ \advance\dimen@\dp\z@
    \advance\dimen@\insertskipamount
    \advance\dimen@\pagetotal\advance\dimen@-\pageshrink
    \ifdim\dimen@>\pagegoal\@midfalse\p@gefalse\fi\fi
  \if@mid%
    \ifdim\lastskip<\insertskipamount\removelastskip\insertskip\fi
    \nointerlineskip\box\z@\penalty-200\insertskip
  \else%
    \SaveLastSkip
    \insert\topins{\penalty100 
    \splittopskip\z@skip
    \splitmaxdepth\maxdimen \floatingpenalty\z@
    \ifp@ge \dimen@\dp\z@
    \vbox to\vsize{\unvbox\z@\kern-\dimen@}
    \else \box\z@\nobreak\insertskip\fi}
    \RestoreLastSkip
   \fi\endgroup}


  \newcount\notenumber
  
  \def\note{\advance\notenumber by 1
    \footnote{\the\notenumber)}}

  \newbox\footbox

  \def\footnote#1{\let\@sf\empty
    \ifhmode\edef\@sf{\spacefactor\the\spacefactor}\/\fi
    \sam${}^{\fam0 #1}$\@sf\vfootnote{#1}}%

  \def\vfootnote#1{\insert\footins\bgroup
     \interlinepenalty100 \splittopskip=1pt
     \floatingpenalty=20000
     \leftskip=0pt\rightskip=0pt%
     \parindent=.3em
     \Smallfonts\rm
     \FootItem@{#1}
     \futurelet\next\fo@t}

  \def\FootItem@#1{\par\hangafter1\hangindent=\FootHang
     \setbox0=\hbox{\ignorespaces#1\unskip}%
     \dimen0=.4em\SetOverhang@
     \noindent\rlap{\box0}\kern\Overhang\ignorespaces}


  \def\fo@t{\ifcat\bgroup\noexpand\next \let\next\f@@t
    \else\let\next\f@t\fi \next}
  \def\f@@t{\bgroup\aftergroup\@foot\let\next}
  \def\f@t#1{\baselineskip=10pt\lineskip=1pt
            \lineskiplimit=0pt #1\@foot}%
  \def\@foot{
        \hbox{\vrule height0pt depth5pt width0pt}
        \egroup}
  \skip\footins=12 pt plus 0pt minus 0pt 
  \count\footins=1000 
  \dimen\footins=8in 



 \def\osumess#1{\EdSpider{\immediate\write16{Line \the\inputlineno: #1}}}%
 \def\HideEdStuff{\gdef\EdSpider##1{}}

 \font\BigSym=cmmi10 scaled \magstep 4

 \def\change{\InLMargin{\hbox{\BigSym \char63\kern10pt}}}

 \def\beginchange{\InLMargin{\hbox{\sam\twelvepoint$\heartsuit$\kern10pt}}}

 \def\endchange{\InLMargin{\hbox{\sam\twelvepoint$\spadesuit$\kern10pt}}}

 \def\InLMargin#1{\strut\vadjust{%
     \kern-\strutdepth
     \vtop to \strutdepth{%
         \baselineskip\strutdepth
         \llap{\sam$\smash{\hbox{\EdSpider{#1}}}$}\null}}}

 \def\strutdepth{\dp\strutbox}
 \def\strutheight{\ht\strutbox}

 \def\NoteInRMargin#1{\strut\vadjust{%
     \kern-1.001\strutdepth
     \vtop to \strutdepth{%
       \baselineskip\strutdepth
       \vss\rlap{\ninepoint\unskip\hskip\hsize
         \vtop to 0pt{%
           \hsize=16em\hfuzz=\hsize
           \leftskip=10pt%
           \rightskip=0pt plus 10000pt%
           \baselineskip=9.8pt\lineskip=.2pt%
           \let\\\break
           \noindent\EdSpider{#1}\vss}%
                \kern10pt}\hbox{}}
       }}

 \def\ednote#1{\NoteInRMargin{\tentt #1}}

 \def\cbar{\InLMargin{%
      \dimen0=\strutdepth\advance\dimen0 by \lineskip
      \vrule width 3pt
      height \strutheight depth \dimen0 \kern
      3pt}}

 \def\ccbar{\InLMargin{%
      \dimen0=2\strutdepth\advance\dimen0 by 2\lineskip
      \vrule width 3pt
        height 3\strutheight depth \dimen0 \kern
      3pt}}

 \newinsert\TRMargIns
 \dimen\TRMargIns=\maxdimen

  \def\Ednote#1{\insert\TRMargIns{%
       \vbox to 0pt{\hsize=140pt\hfuzz=\hsize
           \leftskip=6pt%
           \rightskip=0pt plus 10000pt%
           \baselineskip=9.8pt\lineskip=.2pt%
           \let\\\break
           \SetPageRemainder
           \vglue540pt\vglue-\PageRemainder
           \noindent\EdSpider{\tentt #1}\vss}%
       \smallskip}}

 \def\KillEdStuff{\def\ednote##1{}\def\Ednote##1{}%
      \let\change\relax\let\beginchange\relax\let\endchange\relax
       \let\cbar\relax\let\ccbar\relax}


  \topskip=12pt
  \newskip\StdBaselineskip 
  \StdBaselineskip 12pt
  \lineskip=1.1pt
  \lineskiplimit=.8pt
  \widowpenalty=10000 
  \clubpenalty=10000  
  \abovedisplayskip=6pt plus 1pt minus 1pt
  \abovedisplayshortskip=3pt plus 1.5pt
  \belowdisplayskip=6pt plus 1pt minus 1pt
  \belowdisplayshortskip=5pt plus 1pt minus 1pt
  \hfuzz=1.5pt   

  \def\StdPretolerance{100}
  \tolerance=\StdPretolerance

  \newdimen\StdMathsurround
  \StdMathsurround=1.5pt 
  \mathsurround=\StdMathsurround
  \Mas                   

   \def\prose{\relax\hbox{\kern.6\StdMathsurround}}
  
  \def\StdParskip{0pt}    
  \parskip=\StdParskip
  \parindent=0.5cm
 

  \def\Times{ptmr  } 
  \def\TimesI{ptmri  } 
  \def\TimesB{ptmb  }
  \def\TimesBI{ptmbi  }
  \def\HelveticaN{phvrrn }

  =\Times at 10bp
  =\TimesB at 10bp
  \font\tenit=\TimesI at 10bp
  =\TimesBI at 10bp

  \font\tenmrm=cmr10  


    =\Times at 9bp 
    \font\nineit=\TimesI at 9bp 
    =\TimesB at 9bp 
    =\TimesBI at 9bp 

    =\HelveticaN at 9bp 


  =\Times at 12bp
  \font\twelveit=\TimesI at 12bp
  =\TimesB at 12bp


  \font\titleit=\TimesI at 14.4bp
  =\TimesB at 14.4bp

 \SetAuthorHead{AuthorHead} 
 \SetTitleHead{TitleHead}  


  \def\lBr{\raise.125ex\hbox{[\kern.1125ex}}
  \def\rBr{\raise.125ex\hbox{\kern.1125ex]}}

 \setbox\footbox=\hbox{\Smallfonts 2)~}



  \bgroup
  \catcode`\@=11 
  \gdef\itSpacing{%
     \xspaceskip=.31em plus.1em minus.05em \sfcode `f=2001
     \itWarning@\let\itWarning@\itWarning@@}
  \gdef\itSpacingOff{%
     \xspaceskip=0pt \sfcode `f=1000
     \let\itWarning@\relax}
   \global\let\itWarning@\relax
  \gdef\itWarning@@{\errmessage{%
  Special italic spacing already in force
  (you have probably omitted an ``endth'').
  See itSpacing macro in osuPSfnt.sty
         }}
  \egroup

 \fontdimen1\titlebf=0.0pt
 \fontdimen2\titlebf=3.6135pt
 \fontdimen3\titlebf=2.8908pt
 \fontdimen4\titlebf=1.44539pt
 \fontdimen5\titlebf=6.64882pt
 \fontdimen6\titlebf=14.45398pt
 \fontdimen7\titlebf=1.60439pt

 \fontdimen1\tenbi=0.26794pt
 \fontdimen2\tenbi=2.50937pt
 \fontdimen3\tenbi=2.00749pt
 \fontdimen4\tenbi=1.00374pt
 \fontdimen5\tenbi=4.59717pt
 \fontdimen6\tenbi=10.03749pt
 \fontdimen7\tenbi=1.11415pt

 \fontdimen1\twelverm=0.0pt
 \fontdimen2\twelverm=3.01125pt
 \fontdimen3\twelverm=2.409pt
 \fontdimen4\twelverm=1.2045pt
 \fontdimen5\twelverm=5.39615pt
 \fontdimen6\twelverm=12.045pt
 \fontdimen7\twelverm=1.33699pt

 \fontdimen1\twelveit=0.27731pt
 \fontdimen2\twelveit=3.01125pt
 \fontdimen3\twelveit=2.409pt
 \fontdimen4\twelveit=1.2045pt
 \fontdimen5\twelveit=5.37207pt
 \fontdimen6\twelveit=12.045pt
 \fontdimen7\twelveit=1.33699pt

 \fontdimen1\twelvebf=0.0pt
 \fontdimen2\twelvebf=3.01125pt
 \fontdimen3\twelvebf=2.409pt
 \fontdimen4\twelvebf=1.2045pt
 \fontdimen5\twelvebf=5.5407pt
 \fontdimen6\twelvebf=12.045pt
 \fontdimen7\twelvebf=1.33699pt

 \fontdimen1\tenrm=0.0pt
 \fontdimen2\tenrm=2.50937pt
 \fontdimen3\tenrm=2.00749pt
 \fontdimen4\tenrm=1.00374pt
 \fontdimen5\tenrm=4.49678pt
 \fontdimen6\tenrm=10.03749pt
 \fontdimen7\tenrm=1.11415pt

 \fontdimen1\tenit=0.27731pt
 \fontdimen2\tenit=2.50937pt
 \fontdimen3\tenit=2.00749pt
 \fontdimen4\tenit=1.00374pt
 \fontdimen5\tenit=4.47672pt
 \fontdimen6\tenit=10.03749pt
 \fontdimen7\tenit=1.11415pt

 \fontdimen1\tenbf=0.0pt
 \fontdimen2\tenbf=2.50937pt
 \fontdimen3\tenbf=2.00749pt
 \fontdimen4\tenbf=1.00374pt
 \fontdimen5\tenbf=4.61723pt
 \fontdimen6\tenbf=10.03749pt
 \fontdimen7\tenbf=1.11415pt

 \fontdimen1\ninerm=0.0pt
 \fontdimen2\ninerm=2.25842pt
 \fontdimen3\ninerm=1.80673pt
 \fontdimen4\ninerm=0.90337pt
 \fontdimen5\ninerm=4.0471pt
 \fontdimen6\ninerm=9.03374pt
 \fontdimen7\ninerm=1.00273pt

 \fontdimen1\nineit=0.27731pt
 \fontdimen2\nineit=2.25842pt
 \fontdimen3\nineit=1.80673pt
 \fontdimen4\nineit=0.90337pt
 \fontdimen5\nineit=4.02904pt
 \fontdimen6\nineit=9.03374pt
 \fontdimen7\nineit=1.00273pt

 \fontdimen1\ninebf=0.0pt
 \fontdimen2\ninebf=2.25842pt
 \fontdimen3\ninebf=1.80673pt
 \fontdimen4\ninebf=0.90337pt
 \fontdimen5\ninebf=4.15552pt
 \fontdimen6\ninebf=9.03374pt
 \fontdimen7\ninebf=1.00273pt


 \newcount\MaxSpaceFactor
 \MaxSpaceFactor=3000 

 \def\ItemStyle{\rm}
 \def\NrStyle{\rm}
 \def\ItemItemStyle{\rm}

 \MaxItemTag{(iii)}
 \MaxItemItemTag{(iii)}
 \MaxNrTag{(2)}
 \MaxFootTag{2)}
 \def\ReferenceHang{30pt}

 \catcode`\@=\active


\loadbold

=\Times  
=\Times scaled750
=\Times scaled650
\font\rms=\Times scaled 920 

=\TimesBI scaled 860
=\TimesI scaled 860

\textfont0=\rrm  
\scriptfont0=\erm 
\scriptscriptfont0=\srm

\def\Augment#1#2{%
    \toks0\expandafter{#1}\toks2{#2}%
    \edef#1{\the\toks0\the\toks2}}

 \font\twelverma=\Times  scaled 1200
 \font\tenrma=\Times  scaled 1000
 \font\ninerma=\Times scaled 920
 =\Times scaled 840
 \font\sevenrma=\Times scaled 760
 =\Times scaled 680
 \font\fiverma=\Times scaled 600

 \Augment\tenpoint{%
  \textfont0=\tenrma  \scriptfont0=\sevenrma  
  \scriptscriptfont0=\fiverma  }

 \Augment\ninepoint{%
  \textfont0=\ninerma  \scriptfont0=\sevenrma 
  \scriptscriptfont0=\fiverma}

 \Augment\twelvepoint{%
  \textfont0=\twelverma  \scriptfont0=\ninerma  
  \scriptscriptfont0=\sevenrma}

\mathsurround=1pt
\hsize=13.45truecm
\vsize=19.5truecm
\hoffset=1.25truecm
\voffset=2truecm
\advance\baselineskip by 2pt

\predefine\til{\~}
\def\~#1{\relax\ifmmode\widetilde{#1}\else\til{#1}\fi}

\redefine \le{\leqslant}
\redefine \ge{\geqslant}
\define \wt#1{\mathaccent"0365{#1}}
\define \wh#1{\mathaccent"0362{#1}}

\define \iss{\,\Mathaccent{\raise -.8 ex\hbox{$\widetilde{}$\kern.1em}}\rightarrow\,}

\define \prlim{{\varprojlim}\vphantom{i}\,}
\define \inlim{{\varinjlim}\vphantom{i}\,}

\define\Car{\mathop{\fam0 C}}

\define \ur{\mathop{\fam0 ur}}

\define \id{\operatorname{\fam0 id\,}}
\define \im{\mathop{\fam0 im}}
\define \tpp{\mathop{\fam0 top}}
\define \ab{\mathop{\fam0 ab}}

\define \sep{\mathop{\fam0 sep}}

\define \coker{\mathop{\fam0 coker}}

\define \chr{\mathop{\fam0 char}\,}

\define \Br{\operatorname{\fam0 Br}}

\define \Gal{\mathop{\fam0 Gal}}
\define \Hom{\operatorname{\fam0 Hom}}

\def\bitem{\item{$\bullet$}}

\Mas
\HideEdStuff
\rm 
 

\def\issn{{\nineit ISSN 1464-8997 (on line) 1464-8989 (printed)}}

\def\gtp{{\nineit Published 10 December 2000: \ \copyright\ Geometry \& 
Topology Publications}}

\def\gtv3{{\nineit Geometry \& Topology Monographs, Volume 3 (2000) --
Invitation to higher local fields}}


\def\lione
{{\rms Geometry \& Topology Monographs}}

\def \litwo{{\rms Volume 3: Invitation to higher local fields
}} 

\def\tinfo #1.#2.#3-#4
{{
\noindent  {\lione} \hfill 
\par 
\vskip-1.5pt
\noindent {\litwo} \hfill
\par 
\vskip-1,5pt
\noindent {\rms Part #1, section #2, pages #3--#4} \hfill
\vskip24pt 
}}

\def\tinfos #1.#2.#3-#4
{{
\noindent  {\lione} \hfill 
\par 
\vskip-1.5pt
\noindent {\litwo} \hfill
\par 
\vskip-1.5pt
\noindent {\rms Pages #3--#4} \hfill
\vskip24pt 
}}

\def\tinfoi #1
{{
\noindent  {\lione} \hfill 
\par 
\vskip-1.5pt
\noindent {\litwo} \hfill
\par 
\vskip-1.5pt
\noindent {\rms Pages iii--xi: Introduction and contents} \hfill
\vskip26pt 
}}


  \def\titlepagehead{\hfil}

  \newif\iftitlepage\titlepagefalse
  \newif\ifblankpage\blankpagefalse
  \def\makeheadline{
     \ifblankpage{}\else%
     \iftitlepage
\vbox{\line{\vbox to 8.5pt{}
\ninerm
\copy\HLinebox \hfill
\hglue5mm\ninebf\folio 
\titlepagehead}}%
      \else
\vbox{\ifodd\pageno\rightheadline\else\leftheadline\fi}%
      \fi\vskip 12pt\fi}%
     \def\rightheadline{\line{\vbox to 8.5pt{}%
      \ninerm
\copy\TitleBox \hfill
\hglue5mm\ninebf\folio}}%
     \def\leftheadline{\line{\vbox to 8.5pt{}%
        \unskip\ninerm\unskip\ninebf\folio\hglue5mm
 \hfill \copy\AuthorBox
}}

 \footline={\ifblankpage{}\else
\iftitlepage\ninepoint\sam\hfill
\line{\vbox to 8.5pt{}
\copy\TFLinebox
\hfill
\hglue5mm 
}
            \else
\ninepoint\sam\hfill
\line{\vbox to 8.5pt{}
\copy\FLinebox
\hfill 
\hglue5mm
}
\hfil\fi\global\titlepagefalse\fi}

\def\blankpage{{\blankpagetrue\noindent\hbox to 10pt{\hss}\vfill
\pagebreak}}

\tenpoint\rm 
 